\documentclass[12pt]{article}
\usepackage{amsmath}
\usepackage{amsfonts}
\usepackage{amssymb}
\usepackage{verbatim}

\newtheorem{Theo}{Theorem}[section]

\newtheorem{Lemm}[Theo]{Lemma}
\newtheorem{Prop}[Theo]{Proposition}
\newtheorem{ass}{Assumption}

\newtheorem{Cor}[Theo]{Corollary}
\newtheorem{pkt}[Theo]{}
\newcommand{\pf}{{\bf Proof}\ }

\newcommand{\R}{{\mathbb R}}

\newcommand{\Q}{{\mathbb Q}}

\newcommand{\Z}{{\mathbb Z}}

\newcommand{\N}{{\mathbb N}}

\newcommand{\C}{{\mathbb C}}
\newcommand{\F}{{\mathbb F}}

\newcommand{\LL}{{\cal L}}

\newcommand{\A}{{\rm  F}^{\times}} 
\newcommand{\V}{{\rm  V}}
\newcommand{\E}{{\cal E}}

\newcommand{\EC}{{\cal EC}^K}

\newcommand{\D}{{\mathbb F}^K}

\newcommand{\codim}{{\rm codim\ }}

\newcommand{\card}{{\rm card\ }}

\newcommand{\pr}{{\rm pr}}

\newcommand{\ex}{\exp}
\newcommand{\acl}{{\rm acl}}

\newcommand{\ssp}{{\rm span}}

\newcommand{\qftp}{{\rm qftp}}
\newcommand{\stp}{{\rm sttp}}
\newcommand{\ssn}{\section}

\newcommand{\df}{{\bf Definition}\ }

\newcommand{\bl}{\begin{Lemm}}

\newcommand{\el}{\end{Lemm}}

\newcommand{\bt}{\begin{Theo}}

\newcommand{\et}{\end{Theo}}

\newcommand{\bp}{\begin{Prop}}

\newcommand{\ep}{\end{Prop}}
\newcommand{\bpk}{\begin{pkt}\rm }
\newcommand{\epk}{ \end{pkt}}
\newcommand{\ba}{\begin{ass}}
\newcommand{\ea}{\end{ass}}
\newcommand{\be}{\begin{equation}}
\newcommand{\ee}{\end{equation}}

\newcommand{\fr}{^{\smallfrown}}
\newcommand{\bc}{\begin{Cor}}

\newcommand{\ec}{\end{Cor}}

\newcommand{\lb}{\label}

\newcommand{\ra}{\rangle}

\newcommand{\la}{\langle}

\newcommand{\subs}{\subseteq}

\newcommand{\ld}{{\rm lin.dim}}
\newcommand{\trd}{{\rm tr.deg}}
\newcommand{\qed}{$\Box$ \\}
\newcommand{\SCH}{{\rm SCH}}
\newcommand{\smin}{\setminus}

\newcommand{\sups}{\supseteq}

\begin{document}
\title{The theory of exponential sums}
\author{ 
Boris Zilber\\
University of Oxford}

\abstract{Author: B. Zilber\\ 

Title: The theory of exponential sums\\ 



We consider the theory of algebraically closed fields of characteristic zero with multivalued operations $x\mapsto x^r$
(raising to powers). It is in fact the theory of equations in exponential sums.

In an earlier paper we have described  complete first-order theories of such  structures, conditional on a Diophantine conjecture that generalises
the Mordell-Lang conjecture (CIT). 
Here we get this result unconditionally.  

The field of complex numbers with raising to real powers  satisfies the corresponding theory if Schanuel's conjecture holds. In particular, we have proved that  a (weaker) version of 
Schanuel's conjecture
implies that every well-defined system of exponential sums with real exponents has a solution in the complex numbers. Recent result by Bays, Kirby and Wilkie states that the required 
version of Schanuel's conjecture holds for almost every choice of exponents. It follows that for the corresponding choice of real exponents 
we have an unconditional description of the first order theory of the complex numbers with raising to these powers.
 }
\date{20 December, 2011}

\maketitle

\ssn{Introduction} 
An exponential sum  in variables $x_1,\ldots,x_n$ over the complex numbers with exponents in $K\subs \C$ has the form
$$f(x)=\sum_{i=1}^m a_{i}\exp r_{i}x$$
where $a_{i}\in \C$ and $r_{i}x=\sum_{j=1}^n r_{ij}x_j,$ $r_{ij}\in K.$
When $K=\Z$ this is equivalent to a Laurent polynomial  in variables $y_j=\exp x_j.$ 

Systems of  exponential sum equations with $K=\R$ were studied in \cite{Ko},\cite{Ka}, \cite{Kh} by differential-geometric methods. In \cite{Z2} we started a study of the general case using model theoretic approach. In \cite{Z1} we have introduced a system of axioms for an abstract structure based on an algebraically closed field $\F$ of characteristic zero with a formal operation $x\mapsto \exp rx,$ $r\in K,$ in which exponential sums with exponents in $K$ make sense (a field with raising to powers.). 
This  theory given axiomatically  resolves all the algebraic questions about systems of exponential sums, such as when a system has at least $n$ distinct solutions. 
The problem with the axiom system in \cite{Z1} is that it only makes sense under  a diophantine conjecture CIT
(conjecture on intersections of algebraic varieties with tori) stated in that paper, which was also independently formulated in an equivalent form by E.Bombieri, D.Masser and U.Zannier in \cite{BMZ}. This conjecture generalises the Mordell-Lang conjecture and at present remains open, besides some special cases.

The main result of the present paper\footnote{The earlier version of the paper has been in circulation since 2004 under the title ``Raising to powers revisited''.} is a reformulation of the axioms for a field with raising to powers $K$ in such a way that CIT is not needed. Instead we use M.Laurent's theorem \cite{La} that proves Lang's conjecture for the group $\mathbb{G}_m^n(\C),$ along with a theorem by J.Ax \cite{Ax} stating a form of Schanuel's conjecture for differential fields. We prove that  the axiom system is consistent, complete and for each $K$ the theory is model-theoretically quite "tame", or more precisely, the theory is superstable and near model complete. The latter can be reformulated in  geometric terms: Start with the family of {\em projective} subsets of $\F^n,$ all $n,$ obtained by applying  projections $\pr:\F^n\to \F^{n-m}$ to
 zero-sets of exponential sums in $\F^n$ and construct new sets applying Boolean operations to projective sets.
The resulting family of sets will be closed under projections. 
      
What are possible applications of these result outside model-theoretic context? An obvious idea is to try to establish that the axioms hold for $\F=\C,$ at least for some special cases of this. Results in this direction are proved in this paper. Before formulating these we discuss the acclimatization.      

The first observation systematically followed in \cite{Z2} is that the theory of exponential sums heavily depends on Schanuel's conjecture. A corollary to Schanuel's conjecture relevant to raising to powers can be formulated as follows. Let $K\subset \C$ be
a subfield of finite transcendence degree $d=\trd(K).$ Then, for any $x_1,\ldots,x_n\in \C$
$$\ld_K(x_1,\ldots,x_n)+\trd(e^{x_1},\ldots,e^{x_n})-\ld_{\Q}(x_1,\ldots,x_n)+d\ge 0,$$
where $\ld_K(x_1,\ldots,x_n)$ is the $K$-linear dimension, and similarly for $\Q.$

This we call the Schanual conjecture for $\C^K,$ where $\C^K$ stands for the field of complex numbers with raising to powers $K\subs \C.$

We prove  (Theorem~\ref{6.2})   that in case  $K\subs \R,$ if 
the corresponding version of Schanuel's conjecture holds, $\C^K$ is {\em exponentially-algebraically closed}, that is any well-defined system of  exponential sums equations with exponents in $K$ has a solution in $\C.$  
 The conjecture that $\C_{\exp},$ the  field of complex numbers  with $\exp,$ is exponentially-algebraically closed was made in \cite{Z3}.
 This has been studied in \cite{HR}, \cite{DM} and elsewhere including attempts to refute the conjecture.
 Theorem~\ref{6.2} brings hopes that in general exponential-algebraic closedness follows from Schanuel's conjecture. 
 So far we don't know if this is true even for $\C^K$ when $K$ is not a subfield of the reals.

Another important corollary of the main theorem
 (Corollary~\ref{oldt2})  states that under Schanuel's conjecture for $\C^K$  solutions to an overdetermined system of exponential sums equations lie in a finitely many cosets of proper $\Q$-linear subspaces. Moreover, there is a bound on the number of such cosets, uniform in coefficients of the system (but possibly not on the exponents).

 Finitely, we invoke a recent result by M.Bays, J.Kirby and A.Wilkie that implies
 that for "almost any" tuple $\lambda$ in $\C,$ for $K=\Q(\lambda),$ the structure $\C^K$ satisfies the corresponding version of 
Schanuel's conjecture. Thus, the above theorems are applicable to such $\C^K$ unconditionally. In particular, when 
also $\lambda\subs \R,$ we know the complete theory of $\C^K.$
 
 \medskip

Of course, the Mordell-Lang conjecture is confirmed in full generality now
and, as in \cite{Z1}, one can easily replace $\A$ by any semiabelian variety $\mathbf{A}$  and carry
 out the same construction and axiomatisation since also a corresponding analogue of Ax's Theorem and its  corollaries  
 is available. More precisely, one needs the following (weak CIT) to hold for $\mathbf{A}.$ 
\bpk   \lb{ABK} {\bf Theorem} (J.Kirby, \cite{K})
{\em  Let $\mathbf{A}$ be a complex semiabelian variety of dimension $g$ and $\ex: \C^g\to \mathbf{A}$ the universal covering map. 
Given $W(e)\subs \mathbf{A}^{n},$ with $W\subs  \mathbf{A}^{n+\l}$ algebraic subvariety defined  over $k$
and 
$e\in \mathbf{A}^{\l}$ there are finitely many  codimension $1$ \  ${\rm End}\mathbf{A}$-linear subspaces 
$\mu(W)=\{ M_1,\dots, M_m\}$ of $\C^{ng}$  
 such that   for any ${\rm End}\mathbf{A}$-linear subspace $N\subs \C^{ng},$  $b\in
\mathbf{A}^{\l},$ 
  and
any positive dimensional  atypical component $S$ of the intersection
$W(e)\cap \exp(N)\cdot b$ there is $M\in \mu(W)$ and $s\in S$ such that
$S\subs \exp(M)\cdot s.$}

Here, {\em atypical } for an irreducible component $S$ of the intersection of algebraic subvarieties $W(e)$ and $\ex(N)\cdot b$ 
of $\mathbf{A}^n$ (observe that $\ex(N)$ is an algebraic subgroup of $\mathbf{A}^n$) means that
$$\dim S>\dim W(e)+ \dim \exp(N)-\dim \mathbf{A}^n.$$
\epk
 
To prove our main result in full generality, for semi-abelian $\mathbf{A},$ we would need to consider $\C^g$ as an  
${\rm End}\mathbf{A}$-module which, for $g>1,$ contains divisors of zero and makes the linear structure on $\C^g$ 
more involved. 
This would complicate  definitions and some argument, without visible advantages for applications. But the proof below goes through
practically without changes for any elliptic curve without complex multiplication defined over $\Q.$

\ssn{Preliminaries}

\bpk This section along with definitions and notations discusses basic
ingredients of Hrushovski's construction which is standard enough, so
the reader can guess the proofs if they seem too short or are absent.
 
We use here some of the terminology of \cite{Z2}, slightly improved, where
we discussed $K$-linear and affine spaces, tori and their
intersections with algebraic varieties. 
 
For technical reasons we find it more convenient to represent the
two-sorted structures $(\V,{\rm F})$ with maps $\ex: \V\to {\rm F}$  in the
equivalent way as
 one sorted structures in the language $\LL_K$ which is the extension of
the language of vector spaces over $\Q$ by:
\begin{itemize}
\item 
an equivalence relation $E(x,y)$ meaning $\exp x=\exp y,$
\item
$n$-ary predicates $L(x_1,\dots,x_n)$ for linear subspaces $L\subs \V^n$
given by a set of $K$-linear equations in $x_1,\dots,x_n,$

\item
an $n$-ary predicate $\ln W,$  the  pull-back of (the set of ${\rm F}$-points) an 
algebraic variety $W\subs (\A)^n,$ for each $W$ definable  over $\Q.$\end{itemize}

   \epk
   
\bpk \lb{PK}
\df
$\E^K$ is the class of structures $\D$ in language $\LL_K$ with axioms
saying that $\V$ is
an infinite-dimensional vector space over $K,$  \ $E$ is an equivalence relation on $\V,$ 
$\A=\V/E$  can be identified with  the multiplicative group $\A$ of a  
field ${\rm F}$ of characteristic zero,  
 and the predicates $\ln W$ are inverse images of
algebraic varieties over $\Q$ under the canonical mapping
$$\ex: \V\to \A.$$ 
We also postulate that $\exp$ is a surjective homomorphism from the additive group $(\V,+)$ onto the multiplicative group $\A$ of the field, and 
$\ker,$ the kernel  of the homomorphism is a $\Z$-group, that is the additive group which has the same first-order theory as the group $(\Z,+).$\\

We denote ${\rm PK}$  the underlying set of axioms   ({\em powered field with exponents in $K$}). \epk

\bpk {\bf Notations}
For finite $X, X'\subs \V$ and  $Y,Y'\subs \A$ denote

$\ld_K(X), $ the  dimension of the vector space $\ssp_K(X)$
 generated by $X$ over $K;$ \

$\ld_{\Q}(X), $ the  dimension of the vector space $\ssp_{\Q}(X)$ generated by
$X$ over $\Q;$

$\trd(Y),$ the transcendence degree of $Y;$

$\delta^K(X),$ {\bf the predimension} of finite $X\subs \V:$
$$\delta^K(X)=\ld_K(X)+\trd (\ex (X))- \ld_{\Q}(X);$$
$$\delta^K(X/X')=\delta^K(X\cup X') - \delta^K(X');$$ 

More generally for arbitrary $X,X'\subs \V$ and $k\in \Z,$  there is a usual natural meaning to the expressions $\delta^K(X/Z)\ge k$ and  $\delta^K(X/Z)\ge k.$
 \\

We recall that $A\subset \V$ is said to be {\bf self-sufficient in $\D$} if $\delta^K(X/A)\ge 0$ for 
all finite $X\subs \V.$ This is written as $$A\le \D.$$
We also write $a\le \D$ for tuples $a=\la a_1,\ldots,a_n\ra\in \V^n,$ meaning  
$\{ a_1,\ldots,a_n\}\le \D.$

\medskip

Note that by our assumptions $\V$ is the universe of $\D,$ so we often say $A$ is a subset in $\D$ instead of
$A\subset \V.$

\epk
\bpk
We use 
$\ker$ to denote the kernel of $\exp$ and also 
 the name of the corresponding unary predicate of type $\ln W$  ($W({\rm F})=\{ 1\}$).
 For a substructure $A\subset \D$ we write $\ker_{|A}$ for the realisation of this predicate in $A,$ that is for
 $\ker\cap A.$\\ 

Given $d\in \Z, $ let $\E_d^K$ be the subclass of $\E^K$ consisting of all 
$\D$
satisfying the condition:
   $$\delta^K(X/\ker)\ge -d\mbox{ for all finite }X\subs \V,$$
where $\ker=\ker_{|\F}.$   \epk




\bpk {\bf Remark} For any $\D\in \E_d^K$   
and  $X\subs \ker_{|\F},$ 
$$\delta^K(X)= 0$$
and thus $\E^K_d$ is empty for $d<0.$\\ 

On the other hand, for  any $K$ we have by
Lemmas 2.7 and 2.8 of \cite{Z1}
 $$\E^K_0\neq \emptyset.$$\epk

\bpk \label{caseC}
Assuming $\A=\C^{\times},$ the algebraic torus, and  the Schanuel conjecture we can make a better estimates for the 
minimal $d$ such that $\C^K,$ the complex numbers with raising to powers
$K\subs \C,$ belongs to $\E^K_d.$ 

By Schanuel's conjecture, for any finite $X\subs \C$
$$\trd(X,\exp(X))\ge \ld_{\Q}(X).$$
Recall that we assumed that $\trd(K)$ is finite. Obviously, 
$$\ld_K(X)+\trd(\exp X)+\trd(K)\ge \trd(X,\exp(X)).$$
Hence, 
$$\delta^K(X)\ge \trd(X,\exp(X))-\ld_{\Q}(X)-\trd(K)\ge -\trd(K).$$
Since $\ld_K(\ker)=\ld_{\Q}(\ker)=1$ for $\ker=2\pi i\Z,$ we have 
$$\delta^K(X/\ker)\ge -(\trd(K)+1).$$
Thus, under the conjecture,
$$\C^K\in \E^K_{d},\mbox{ for }d=\trd(K)+1.$$

 \epk
\bpk \label{caseC+} 
Given $\D\in \E^K_d\smin \E^K_{d-1}$ one can find a finite $A$ in $\D$ with
\be \label{A}\delta^K(A/\ker)=-d.\ee
By minimality $$A\cup \ker \le \D.$$

The equality (\ref{A}) does not change if we extend $A$ by elements of $\ker.$
In case $\delta^K(X)$ is bounded from below by, say $-d',$   for all finite $X\subset \ker,$ in particular, if $\ker$ is a finite rank group, then 
$$\delta^K(Y)\le -(d+d')\mbox{ for all finite $Y$ in } \D.$$  
It follows that there
is a finite $A_0\subs \ker$ such that if $A_0\subs A$ and (\ref{A}) holds then the value of $\delta^K(A)$ reaches minimum 
 and so \be \label{A_0} A\le \D\ee
 that is $$\delta^K(X/A)=\delta^K(X/A\cup \ker)\ge 0\mbox{ for every $X$ in } \D.$$
\epk
\bpk \label{assA} By \cite{Z1} for every field of powers $K$ there is a $d\ge 0$ and there exists a  $\D\in \E^K_d\smin \E^K_{d-1}$ and so 
 $A$ with the property (\ref{A}) and (\ref{A_0})  does exist. Assuming $\D=\C^K$ and  Schanuel's conjecture holds we can choose $d$ and $A$ as above and   (\ref{A_0}) will hold.
In fact, this is a form of Schanuel's conjecture for $\C^K$ for a given $K\subs \C.$

\epk
\bpk \df A structure $\D$ in $\E^K_d$ is said to be {\bf $\E^K_d$-exponentially-algebraically closed}
(e.a.c.) if for any $\D_1\in \E^K_d,$ such that $\D\le \D_1,$  any  finite
quantifier-free type
over $\D$ which is realized in $\D_1$   has a realization in $\D.$

$\EC_d$ will stand for the class of $\E^K_d$-exponentially-algebraically
closed structures, or, in the shorter form, $\EC.$\epk

Using the standard Fraisse construction in the class $\E^K_d$ relative to the
strong embedding 
$\le$ one proves:

\bpk {\bf Propositon} [Proposition 1 of \cite{Z1}] For any $\D$ in $\E^K_d$ there exists an $\E^K_d$-e.a.c. structure
containing $\D.$\medskip

More difficult is the following result, proved using Ax's Theorem (in a more general form one can use Theorem~\ref{ABK}). 
\bpk \label{2.9}{\bf Proposition}  [Corollaries 1 and 2, section 4 of \cite{Z1}]\lb{c3}  There exists a set {\rm EC} of first order $\forall \exists$ axioms such that,  for any $\D\in \E^K_d,$
$$\D\models {\rm EC}\mbox{ iff $\D$ is exponentially-algebraically closed.}$$

\epk
\ssn{The linear structure on \V}

\bpk \lb{not} 

A subspace $L\subs \V^n$ is said to be $K$-{\bf linear}
if there are $ k_{ij}\in K$ \ $(i\le r, \ j\le n)$ \ such that
$$L= \{ \la x_1,\dots,x_n\ra \in \V^n: k_{i1}x_1+\dots+k_{in}x_n=0 \}.$$
Define $\dim L=\mbox{co-rank}(k_{ij}),$ 
the co-rank of the matrix $(k_{ij}).$ \\

Let $L\subs \V^{n+l}$ be a $K$-linear subspace, $a=\la a_1,\dots,a_{\l}\ra.$
Let
$$L(a)=\{ \la x_1,\dots,x_{n}\ra \in \V^n:  \la
x_1,\dots,x_{n},a_1,\dots,a_{\l}\ra\in L\} .$$

We call such an $L(a)$ 
a $K$-{\bf affine subspace defined
over } $a.$ \\ 
The same terminology is applied for $\Q$ instead of $K.$

\epk

\bpk \lb{linalg} {\bf Lemma}  A $K$-affine subspace $L(a)\subs \V^n$ can be
represented equivalently and uniformly on $a$ as
$$L(a)=L(0)+r(a), \ \ r(a)\in  \V^n, \ \ r\mbox{ is a $K$-linear map},\ \ \V^{\l}\to \V^n$$ and $L(0)$ is
a $K$-linear subspace.

Moreover, if  $r'$ is any  $K$-linear map such that $r'(a)\in L(a)$ for all
 $a\in \pr_{n+1\dots n+l} L,$ then also $$L(a)=L(0)+r'(a).$$

\pf
Let $L$ be determined by the system
of linear equations 
\be \lb{N1}\sum_{j=1}^n q_{ij}v_j+\sum_{s=1}^lk_{is}w_{s}=0, \ \
 q_{ij}, k_{is}\in K.\ee 
Then the system
$$\sum_{j=1}^n q_{ij}v_j+\sum_{s=1}^lk_{is}a_{s}=0, $$
determines $L(a).$ It follows that
\be \lb{Nh}\sum_{j=1}^n q_{ij}v_j=0\ee
determines $L(0).$

By linear algebra there is an $n$-tuple $\{ r_j(w_1,\dots,w_l): j=1,\dots, n\}$
of $K$-linear functions $V^l\to V$
 such that, if for a given 
$(w_1,\dots,w_l)$ the system (\ref{N1}) is consistent, then
$v_j=r_j(w_1,\dots,w_l)$ gives a solution to the system. Hence,
$v-r(a)$ is a solution to the homogeneous system (\ref{Nh}) iff $v$ is a solution of the system (\ref{N1}) with $w=a.$

The 'moreover' statement follows immediately from the fact that $r'(a)-r(a)\in L(0).$
\qed
\epk

\bpk \lb{N} {\bf Lemma}  Given a $K$-linear $L\subs V^{n+\l}$ and $0$ the zero of $V^{\l} $
\begin{itemize}
 \item[(i)] there exists a unique maximal $\Q$-linear subspace $N_L\subs L;$
\item[(ii)] $N_{L(0)}=N_L(0);$
\item[(iii)] given
$a\in \pr_{n+1\dots n+\l} L$ and $q(a)  \in L(a)\cap \ssp_Q(a)$ 
 there exists a 
maximal $\Q$-affine 
subspace   $N_{L,q(a)}(a)\subs L(a)$ over $\ssp_Q(a)$ containing $q(a),$ 
 and in this case  $N_{L,q(a)}(a)=N_L(0)+q(a).$ 
\end{itemize}
\medskip

\pf (i) 
$N_L(0)$ exists since the sum of two  $\Q$-linear subspaces of 
$L(0)$ is again a  $\Q$-linear subspace of $L(0).$ 

(ii) Obviously, $N_L(0)$ is a $\Q$-linear subspace of 
$L(0),$ so $N_{L(0)}\sups N_L(0).$ 

 $N_{L(0)}$ is a $\Q$-linear subspaces of $L(0),$ so $N_{L(0)}\times \{ 0\}$
is  a $\Q$-linear subspaces of $L,$ hence $N_{L(0)}\times \{ 0\}\subs N_L$ and
$N_{L(0)}\subs N_L(0).$

(iii)  
 $N_L(0)+q(a)$ 
 is a $\Q$-affine subspace of $L(0)+q(a)=L(a).$   
If $M+q(a)$ is another $\Q$-affine subspace of $L(a),$ containing $q(a)$ then
  $M=(M+q(a))-q(a)\subs L(0)$ and hence $M\subs N_L(0),$    
$M+q(a)\subs N_L(0)+q(a).$
\qed
\epk


\bpk {\bf Remark}
Note that given a $\Q$-linear subspace $N$ of $\V^n,$ its image $\ex(N)$ is an algebraic subgroup
of $(\A)^n,$ a torus. Correspondingly, the image of a $\Q$-affine subspace is a coset of 
an algebraic subgroup and so has a form $\ex(N)\cdot b$ for some $b\in (\A)^n.$

\epk

\ssn{Intersections with cosets of algebraic subgroups}

For $W\subs (\A)^{n+\l}$  an algebraic variety, $b=\la b_1,\dots,b_{\l}\ra$
denote
$$W(b)=\{ \la x_1,\dots,x_{n}\ra \in (\A)^n:  \la
x_1,\dots,x_{n},b_1,\dots,b_{\l}\ra\in W \}.$$

Below we  refer to Theorem~\ref{ABK} but will only use the case when the semi-abelian variety in question
is just the multiplicative group $\A$ of the field,  which is a theorem in \cite{Z2} and in an equivalent form in \cite{BMZ}.

\bpk \lb{fact0} {\bf Lemma.} In the statement of  Theorem~\ref{ABK} we can assume that
$s\in \acl(b,e).$
\smallskip

\pf $S$ is an irreducible component of the set $W(e)\cap \ex(N)\cdot b$ 
definable over $(e,b)$ hence it is definable over $\acl(b,e).$ Thus it
contains points from the algebraically closed field $\acl(b,e).$\qed
\epk

\bpk \lb{fact} {\bf Proposition.}
Given $W\subs (\A)^{n+\l},$ an algebraic subvariety defined  over $\Q$
 there are finitely many proper $\Q$-linear subspaces 
$\pi(W)=\{ M_1,\dots, M_m\}$ of $V^n$  
 such that   for any $e,b\in (\A)^l$ and a 
$\Q$-linear subspaces   $N\subs V^{n},$ for any  
positive dimensional  atypical component $S$ of the intersection
$W(e)\cap \ex(N)\cdot b$ there  is a
$M\in \pi(W)$ and $s\in S\cap \acl(e,b)$ such that $S\subs \ex(M)\cdot s$ and 
 $S$  is a typical component of 
$\ex(N)\cdot b\cap \ex(M)\cdot s\cap W(e)$ with respect to the group variety $\ex(M)\cdot s.$
\smallskip

\pf Notice first that by obvious transformations of $W$
we can assume that the family $\{ W(e): e\in (\A)^{\l}\}$ is invariant with
 respect to shifts by elements of $(\A)^{\l},$ that is, 
for every $b,e\in (\A)^l$ there is $e'\in (\A)^{\l}$ such that
$$W(e)\cdot b=W(e').$$

By induction on the dimension of a proper algebraic subgroup $P$
of $(\A)^n,$ for any algebraic subvariety $W_P$ of $P\subs (\A)^l$ over
 $\Q$
   we construct a collection of proper 
algebraic subgroups 
$\pi_P(W_P)$ of $P$ such that the statement of the lemma holds for 
$\ex(N)\cdot b\cap W_P(e)\subs P.$

For $\dim P=1$ the statement is trivially true for there is no
 atypical components in any intersection.
 
Consider the general case. Assume by induction that $\pi_P(W_P)$
has been constructed for $\dim P<\dim (\A)^n.$ Notice that by invariance
$\pi_P(W_P)$ will be the same if we replace $P$ by $P\cdot b,$ for
$b\in (\A)^l$ a parameter.

 Assume that for all $P\subset (\A)^n$
proper, $\pi_P(W_P)$ exists.

Given $W\subs (\A)^{n+l},$ 
we 
let
$$\pi(W)=\bigcup_{Q=\ex(M), \ M\in \mu(W)} \pi_Q(W_Q)\cup\{ M\},$$
where $W_Q$ is $W\cap (Q\times (\A)^l).$

Now, if $S\subs \ex(N)\cdot b\cap W(e)$ is atypical, then by Theorem~\ref{ABK}
$S\subs \ex(M)\cdot s$ for some $M\in \mu(W).$ Hence 
$S\subs \ex(N)\cdot b\cap W_Q(e),$ for $Q=\ex(M)\cdot s,$ and is a
component of the intersection. In other words, $S$ is a component of
the intersection $\ex(N\cap M)\cdot b'\cap W_Q(e)$ for some $b'\in (\A)^l.$ Either $S$ is
typical in this intersection with respect to $Q=\ex(M),$ and hence the statement of the
proposition holds for the chosen $M$ belonging to $\pi_Q(W_Q)$ by
definitions, or $S$ is atypical but we can find by induction
$M'\in \pi_Q(W_Q)\subs \pi(W)$ such that $S\subs \ex(M')\cdot s$ and  is a typical component
in the intersection with respect to $\ex(M').$  
  
\qed

Now we want to show  that under certain conditions we can factor out 
$M$ in the previous proposition.\epk

\bpk \lb{3.3}
Let $M\subs \V^n$ be a $\Q$-linear subspace. 
We see $\V^n$  as a subspace of $\V^{n+\l},$ equivalently
$\V^{n+\l}=\V^n\dot+ \V^{\l},$ with $a\in \V^{\l},$ $e=\ex(a)\in (\A)^{\l}.$

By definitions  $M$ is definable by 
$c=\codim M$ independent $\Q$-linear equations 
$$m_{i1}v_1+\dots+m_{in}v_n=0, \ \ i=1,\dots c,$$  
where $(v_1,\dots, v_n)\in \V^{n}.$ The same in matrix  notation
$$\bar m\bar v=\bar 0.$$

We now choose ${\bar m}^{\perp},$ a   $(n-c)\times n$-matrix consisting
of vectors  $(m_{j1},\dots,m_{jn})\in \Q^n,$ $j=c+1,\dots, n,$ which extend $\bar m$ to the basis of the $\Q$-vector space $\Q^n.$
We let $M^{\perp}$ to be the set of solutions to the system 
$${\bar m}^{\perp}\bar v=\bar 0.$$  

This determines  the definable decomposition
$$\V^n=M\dot+M^{\perp}\cong M\dot \times \V^n/M.$$

 Applying $\ex$ we correspondingly have the decomposition
$$(\A)^n= Q\cdot Q^{\perp}\cong Q\times  (\A)^{n}/Q,$$
where $Q=\ex(M)$ and $Q^{\perp}=\ex(M^{\perp}).$ 

Note that the structure
$(M^{\perp},\ex, Q^{\perp})$ is by construction isomorphic to  
 $(V^c,\ex,(\A)^c)$ in the language $\LL_K.$

We denote  the natural mappings $$\V^n\to M^{\perp}\mbox{ and }(\A)^n\to  Q^{\perp}$$
associated with the above decomposition as 
$$v\mapsto v+M\mbox{ and }x\mapsto x\cdot Q,$$
correspondingly. 

It is easy to see that $x\mapsto x\cdot Q$ is a proper mapping on $(\A)^n$
(and $(\A)^{n+l}$)
hence the images $W/\ex M$ (that is $W/Q$) and  
$W(e)/\ex M$ are algebraic subvarieties of $\ex M^{\perp}\times(\A)^{\l} $ and 
$\ex M^{\perp},$ correspondingly. 

The same algebraicity statement holds for the quotients $S/\ex M$ and $b\cdot\ex(N)/\ex M.$  

\epk
\bpk \lb{3.4}
On the other hand, $\ex$ can be naturally extended to the quotient spaces  $$\ex:\ \V^n/M\to (\A)^n/Q.$$ 
So,  $(\V^n/M, \ex, (\A)^n/Q)$ and $(\V^{n+\l}/M, \ex, (\A)^{n+\ell}/Q)$ are canonically isomorphic to 
$(\V^c,\ex,(\A)^c)$ and $(\V^{c+\l},\ex,(\A)^{c+\l})$ correspondingly and hence, for 
$(\V,\ex,\A)\in \E^K_d,$ $u\in \V^n$ and $a\in \V^{\l}, $ 
$a\le \D,$ 
$$\delta^K(u+M)=\ld_K(u+M)+\trd(\ex(u+M))-\ld_{\Q}(u+M)\ge d,$$
and 
$$\delta^K(u+M/a)=\ld_K(u+M/a)+\trd(\ex(u+M)/\ex a)-\ld_{\Q}(u+M/a)\ge 0.$$

\epk
\bpk \lb{3.5} 

We have also the decomposition \be \lb{dec}L(a)=L(\bar 0)\cap M\dot+L/M(a)\ee where 
$$ L/M=L/(L\cap (M\times  \{ \bar 0\}))\subs V^{n+\ell}/(M\times \{ \bar 0\}), \
\bar 0\in \V^{\ell},\ L/M(a)\subs M^{\perp}.$$
Thus we can naturally identify $L(a)/M$ with $L/M(a).$\\

\epk

\ssn{Axiomatizing $\E^K_d.$}   \lb{s4}

Fix $a\in \V^{\l}$ and 
consider the pairs $(L(a),W(\ex a)),$ where $L$ is a $K$-affine subspace of $\V^n$
over $a$
and $W$ an algebraic subvariety of $(\A)^n$ over $\Q.$ \\

\df A pair $(L(a),W(\ex\ a)$ is said to be {\bf special} if $L$ is not contained in any
proper $\Q$-linear subspace of $V^{n+\ell}$  and 
\be \lb{spec} \dim L(a)+\dim W(\ex\ a) < n.\ee

(This corresponds to {\em $0$-special} in the terminology of \cite{Z1}.)

\bpk \lb{d}
 Let, for a $\Q$-subspace $M$ of $\V^n$
  $$d({W(\ex\ a),\ex(M)})=\min \{ \dim  (\bar W(\ex\ a)\cap w\cdot 
\overline{\ex M}: \ w\in
W(\ex a)\},$$
where $\bar W$ and $\overline{\ex M}$ is the closure in the ambient projective space.

 Let 
$$W^{\ex M}(\ex\ a)=\{ w\in W(\ex\ a): \ \dim  (W(\ex\ a)\cap 
w\cdot\ex M)>d({W(\ex\ a),\ex M})\}.$$ 
Since the fibres of minimal dimension  are located over an open subset
 $W^{\ex M}(\ex\ a)$ is a proper closed subset of 
$W(\ex\ a),$ 
  maybe empty, if  $d({W(\ex\ a ),\ex M})=\dim (W(\ex\ a)\cap \ex(x)\ex M).$

\epk
\bpk \lb{w/M}

Suppose that  $(L(a),W\ex\ a))$ is  special. 
Suppose $M\subs L(0).$ Consider the quotients $M^{\perp}=\V^n/M,$\ \
$\ex(M)^{\perp}=(\A)^n/\ex(M)$ and subsets $L(a)/M$ and  
$W(\ex a)/\ex(M).$  

Then,
$\dim W(\ex a)/\ex(M)< n-\dim M,$
that is\\

{\em $W(\ex,a)/\ex(M)$ is a proper subvariety of   $(\A)^n/\ex(M).$}\\

Indeed, by addition formula,
 $$\dim W(\ex a)/\ex(M)= \dim W(\ex a)-d({W(\ex\ a ),\ex M})$$ and 
$$d({W(\ex\ a ),\ex M})\ge 0> \dim W(\ex a)+\dim M-n,$$ since  
$W(\ex a), L(a)$ is special.

\epk

 
\bpk \label{4.3}
{\em Assume $a\le \D$ is and
 let $x\in \V^n.$ 
We  analyse first order consequences of this assumption. We also aim to show  that the analysis yields the same conclusions and formulas when we replace $a$ by $b$ satisfying
the same quantifier-free type,
 $\qftp(a)=\qftp(b).$}\epk \lb{cases}

Suppose 
$(L(a),W\ex\ a))$ is {\bf special}, 
Suppose $x\in L(a)$ and $\ex x\in W(\ex a ),$ in $\D$ and $a\le \D.$

Since
$$\ld_K(x/a)+\trd (\ex(x)/\ex a)-n\le \dim L(a)+\dim W(\ex\ a)-n<0$$
and, as $a\le \D,$ $$\delta^K(x/a)\ge 0,$$
we have $\ld_{\Q}(x/a)<n,$ so
$x\fr a\in N$ for some  proper $\Q$-linear subspace $N$ of $\V^{n+\l}.$ 
We assume $N$ is
minimal for $x\fr a.$

We have $\ex x\in S_x\subs  \ex(N(a))\cap W(\ex a),$ 
where $S_x$ is
a component  of \ $\ex(N(a))\cap W(\ex\ a).$

{\bf Case 1.} The component $S_x$ 
is of dimension $0.$ 

Then 
$\trd(\ex(x)/\ex(A))=0,$ which implies $\ld_K(x/a)=\ld_{\Q}(x/a),$
that is $\dim N(a)=\dim N(a)\cap L(a)$ and so
$N(a)\subs L(a)$ is a $\Q$-affine subspace over $a,$ thus  $N(a)=N_L(0)+q(a),$
 for some
  $q(a)\in L(a)\cap \ssp_Q(a)$ (Lemma~\ref{N}). 

So $$x\in  N_{L}(0)+q(a)\mbox{ for some }
  q(a)\in L(a)\cap \ssp_Q(a).$$
 \\

{\bf Subcase 1.1} $d({W(\ex a),\ex(N_L(0))})<\dim (W(\ex a)\cap  \ex(x+N_L(0))).$

Under this assumption 
$W^{\ex(N_L(0))}(\ex a)$ is a proper closed subset of 
$W(\ex a)$ containing $\ex x,$ by \ref{d}.\\

Otherwise we have 

{\bf Subcase 1.2.} $d({W(\ex a ),\ex(N_L(0))})=\dim (W(\ex a)\cap \ex(x+N_L(0))).$ 

We have $W^{\ex(N_L(0))}(\ex a )=\emptyset$ in this case.

Consider the quotients $N_L^{\perp}(0)=\V^n/N_L(0),$ \
$\ex(N_L(0))^{\perp}=(\A)^n/\ex(N_L(0))$ and subsets $L/N_L(a)$ and  
$W(\ex,a)/\ex(N_L(0)).$  By \ref{w/M} $W(\ex a)/\ex(N_L(0))\subsetneq(\A)^n/\ex(N_L(0)).$

Obviously, for the $x$ above,
$\ex(x/N_L(0))$ is a point in  $W(\ex\ a)/\ex(N_L(0))$  equal to 
$\ex(q(a)/N_L(0)).$ 
\epk

Let
$$\Gamma_a=\{ \ex(q(a)\fr a): q(a)\in \ssp_Q(a)\}$$
This is a coset $\Gamma^0\cdot s(\ex a)$ of a finite rank subgroup 
$\Gamma^0$ of $(\A)^{n+l}$ (depending on the choice of   $\ex a$). 

By Laurent's Theorem  
there  are finitely many, say $k_a,$
cosets $T_i(\ex a)\subs (\A)^{n+\l}$  of group subvarieties (tori) $T_i\subs (\A)^{n+\l},$ 
$T_i(\ex a)\subs W(\ex a)$ (note the notation, $T_i(\ex a)$ as defined in \ref{not})
 such that 
\be \lb{Gamma}\Gamma_a\cap W(\ex a)=\cup_{i\le k_a}\Gamma_a\cap T_i(\ex a).
\ee

Hence, in this case 
\be \lb{c1.2}\ex(x/N_L(0))\in \bigcup_{i\le k_a}T_i(\ex a) \ee
and $T_i(\ex a)\subs W(\ex a)/\ex N_L(0)\subsetneq(\A)^n/\ex(N_L(0).$  \\

{\bf Remark 1.3.} Note that since the ingredients of (\ref{Gamma}) and (\ref{c1.2}) are defined in terms of $\ex\frac{a}{n},$ $n\in \N,$ both (\ref{Gamma}) and (\ref{c1.2})
continue to hold with the same $T_i$ and $k_a$ if we replace $a$ by $b$ with
$\ex\frac{a}{n}\equiv \ex\frac{b}{n},$ all $n,$ in the field language (that is Galois conjugated for each $n.$)\\

{\bf Case 2.} $\dim S_x>0.$ 
Then, by \ref{fact}, $S_x\subs c\cdot\ex(M),$ for some $\Q$-linear subspace 
$M\in \pi_W$ of $\V^{n+\ell},$ $c\in \acl(\ex a),$ and $S_x$ is
typical in the intersection\\ $W(\ex a)\cap c\cdot\ex(M)$ with respect to $c\cdot
\ex(M).$ The latter gives us, for  $a'=\ln c,$ 
$$ \dim S_x=\dim (\ex(N(a)) \cap  \ex(M+a'))+\dim 
(W(\ex a) \cap \ex(M+a'))-\dim \ex(M+a').$$
It is easy to see that $\delta^K(x/aa')\ge 0$ and hence we obtain
$$ \dim (L(a)\cap N(a)\cap (M+a'))+\dim S_x-\dim (N(a)\cap (M+a'))\ge 0.$$

Combining with the above we get
$$\dim (L(a) \cap N(a)\cap (M+a'))+\dim (W(\ex  a)\cap \ex(M+a'))-\dim 
\ex(M+a')\ge 0.$$
And so
\be\lb{typ3}\dim (L(0)\cap M)+\dim (W(\ex  a)\cap \ex(M+a'))-\dim M\ge 0.\ee

{\bf Subcase 2.1} $d({W(\ex  a),\ex(M)})<\dim (W(\ex  a)\cap c\cdot\ex(M)).$  

Under this assumption $W^{\ex(M)}(\ex  a)$ is a proper closed subset of $W(\ex  a)$ containing $\ex(x).$\\

Otherwise, we have

{\bf Subcase 2.2.} $d({W(\ex  a),\ex(M)})=\dim (W(\ex  a)\cap c\cdot\ex(M)).$ 

So, \be \lb{2.2a}W^{\ex(M)}(\ex  a)=\emptyset.\ee 
We now apply the factorisation of \ref{3.3}-\ref{3.4}.
 
Obviously, for our $x,$
 
$$\ex(x/M) \in S_x/\ex(M)\cap \ex(L(a)/M)$$ and $S_x/\ex(M)$ is a singleton in $(\A)^{n-\dim M }$
defined over the same parameters as $S_x,$ that is over $\acl(\ex(A)).$ 

Now we notice that
$$\dim L/M(a)+\dim W(\ex a)/\ex M=$$
$$=\dim L(a)-\dim L(0)\cap M +\dim W(\ex a)-d({W(\ex  a),\ex(M)})=$$
$$=[\dim L(a) +\dim W(\ex a)]-[\dim L(0)\cap M+d({W(\ex  a),\ex(M)})].$$
The sum in the first bracket is less than $n$ by assumptions, and the sum in the second bracket is not less than $\dim M$ by (\ref{typ3}). Hence
$$\dim L/M(a)+\dim W(\ex a)/\ex M<{n-\dim M },$$
that is {\em the  pair is special.}

This means that after factorisation by $M$ we are in case 1 again.

Hence either, as in subcase 1.1,
\be \label{2.2aa}
\begin{array}{ll}
W^{\ex(N_{L/M}(0))}(\ex a)\mbox{ is a proper closed subset of }\\ 
W(\ex a)/\ex M\mbox{ containing }\ex(x+M)\end{array}
                                           \ee
or, as in subcase 1.2,

\be \label{2.2bb}\ex(x+N_{L/M}(0)+M)\in \bigcup_{i\le k_a}T_i(\ex a), \ee
for group subvarieties $T_i,$
$$T_i(\ex a)\subs W(\ex\ a)/\ex (N_{L/M}(0)+M)\subsetneq(\A)^n/\ex(N_{L/M}(0)+M).$$

{\bf Remark 2.3.}
Note again as in  1.3 that (\ref{2.2aa}) and (\ref{2.2bb}) continue to hold with the same $T_i$ if we replace $a$ by $b$ with
$\ex \frac{a}{n}\equiv \ex \frac{b}{n},$ all $n,$ in the field language.

\bpk\label{4.6}
For  $L,$  $W$ and $a$ as in \ref{4.3} let, for  $b$ of the length equal to that of $a,$
$$\begin{array}{lll}\Phi_{L,W,a}(b,x):= x\fr b\in L \  \& \ \ex(x\fr b)\in W \to 
  \bigvee_{M\in \pi_W\vee M=N_L(0) \vee  M=\{ 0\}} \\
\ex x\in W^{\ex M}(\ex b)\ \vee \ \ex(x+M)\in W^{\ex N_{L/M}}(\ex b)\vee \\ 
      \vee \ \ex(x+M)\in \bigcup_{i\le k_a}\ex(N_{i,a,M})  \end{array}$$  
This is a quantifier-free formula without parameters and, by the analysis above 
$$\D\models \forall x \Phi_{L,W,a}(a,x).$$

\epk

\bpk \lb{ingr} {\bf Consequence of the analysis.} Under assumptions \ref{4.3} and \ref{4.6}
we get by Remarks 1.3  and 2.3 (subsection \ref{4.3}):

If $$b\le \D\mbox{ \ and\ }
\ex(\ssp_Qa)\equiv_{\rm fields}\ex(\ssp_Qb)$$
 then 
$$\D\models \forall x \Phi_{L,W,a}(b,x).$$

\epk

\bpk
We define a {\bf strong embedding type} with variables $y=\la y_1,\ldots,y_n\ra,$ where $n$ is the length of $a:$
$$\stp_a(y):= \qftp_a(y)\cup  \{   \forall x \ \Phi_{L,W,a}(y,x):\  \left( L(a),W(\ex a)\right) \mbox{  special} \},$$
where $\qftp_a(y)$ denotes the quantifier-free type of $a$ over $\emptyset$ in variables $Y.$

This is a type consisting of universal formulas without parameters.\epk

Now we can reformulate \ref{ingr} 

\bpk \lb{Claim 1} $b\le \D \mbox{ and } \D\models \qftp_a(b)\Rightarrow  \D\models \stp_a(b).$\epk

{\bf We assume below that $\D \in \E^K_d$ and $a$ has been chosen so that 
$a\le \D$ as well as $a\cup \ker\le \D$} as discussed in \ref{assA}.

\bpk \lb{Claim 3} {\bf Lemma.}

$$  \D\models \stp_a(b) \Rightarrow b\le \D.$$

\pf Assume w.l.o.g. that $x\in V^n$ is $\Q$-linearly
independent over $b,$ 
$L(b)$ is the minimal $K$-affine subspace over $b$ containing $x,$
 and $W(\ex b)$ the minimal algebraic variety over $\ex b$ 
containing $\ex x.$ Notice that under this choice  $\ex x$ is   multiplicatively independent over $\ex b.$

We show that $(L(b),W(\ex  b)$ can not be special, thus proving
$\delta^K(x/b)\ge 0.$ 

Indeed, if the pair were special, 
$\forall x \Phi_{L,W,a}(b,x)$
 would imply that  $x$ satisifes 
 $\ex x\in W^{\ex M}(\ex b)\ \vee \ \ex(x+M)\in W^{\ex N_{L/M}}(\ex b)$ or 
      $\ex(x+M)\in \bigcup_{i\le k_a}\ex(N_{i,a,M})$

Both $\ex x\in W^{\ex M}(\ex b)$ and $\ex(x+M)\in W^{\ex N_{L/M}}(\ex b)$ 
contradict the assumptions that
 $W(\ex b)$ is the algebraic locus of $\ex x$ over $\ex b,$ since
$W^T(\ex b),$ for $T$ a proper algebraic subgroup, is a proper subvariety of 
 $W(\ex b)$ by \ref{d} and (\ref{2.2aa}).

 The condition  $\ex(x+M)\in \bigcup_{i\le k_a}\ex(N_{i,a,M})$ can not hold because
by (\ref{2.2bb}) it would contradict the fact that $\ex x$ is   multiplicatively independent over $\ex b.$
\qed

\epk

 \bpk \lb{p1} {\bf Proposition.} The following two conditions are equivalent:
\be \lb{i} \D\models  \stp_a(b)\ee
and
\be \lb{ii}b\le \D\ \ \&\ \ \D\models \qftp_a(b)\ee

\pf Lemma~\ref{Claim 1} proves (\ref{ii}) $\Rightarrow$ (\ref{i}).
The converse follows from Lemma~\ref{Claim 3} and the definition of $\stp.$\qed
\epk

\bpk Let $\D$ be a member of $\E^K_d$ and $a\le \D$ be a finite tuple containing 
generators of $\ker(\ex)$ such that
$$\delta^K(a)=-d.$$ 
It follows that $a\le \D$  so, we are still under  assumtions of \ref{4.3}.

Let 
$$\SCH_{a}=\{ \exists y \bigwedge S(y) : \ S(y)\subset \stp_a(y)\},$$
be the set of $\exists\forall$-sentences stating the consistency of type $\stp_a(y)$ in variables $y=\la y_1,\ldots,y_n\ra.$

\medskip

Recall \ref{PK} and the notation PK for the axioms of powered fields with powers in K.

\epk

\bpk \lb{Claim 6}{\bf Lemma.}
 Let $\D$ be a model of ${\rm PK +SCH}_{a}$ which realises the type $\stp_a.$ Then $\D\in \E^K_d\smin \E^K_{d-1}.$
\smallskip

\pf
By assumption
 we have  $b$ in $\D$  such that $\D\models \stp_a(b).$
By Proposition~\ref{p1} 
$b\le \D.$ 
As a consequence of $\qftp_a(b)$  we have $\delta^K(b)= -d.$ 
It follows that $\D\notin \E^K_{d-1}.$

To see that $\D\in \E^K_d$ 
 we need to prove that
$\delta^K(Z)\ge -d$  for any finite $Z\subs \D.$

Let $Y$ be a $\Q$-linear basis of
$\ssp_Q(Z)\cap \ssp_Q(b).$  We have then $\ld_{\Q}(Z/b)=\ld_{\Q}(Z/Y)$ and thus
$\delta^K(Z/Y)\ge \delta^K(Z/b)\ge 0.$
But
$\delta^K(Z)=\delta^K(Z/Y)+\delta^K(Y),$ so $\delta^K(Z)\ge \delta^K(Y)\ge -d.$\qed
\epk

\bpk \lb{axio0} {\bf Theorem.}  {\em Assume $\delta^K(a)=-d.$ 
The following two conditions are equivalent for a structure $\D:$
\begin{itemize}
\item[(i)] $\D\models {\rm PK+}\SCH_{a};$
\item[(ii)] $\D\in \E^K_d\smin \E^K_{d-1}$ and $\qftp_a$ is realised in some  ${^*\D}\succ \D.$ 
\end{itemize}

Moreover, also the following two  are equivalent:
\begin{itemize}
\item[(iii)] $\D\models {\rm PK+}\SCH_{a}{\rm + EC};$
\item[(iv)] $\D\in \EC_d\smin \E^K_{d-1}$ and $\qftp_a$ is realised in some  ${^*\D}\succ \D.$ 
\end{itemize}

}
\smallskip

\pf  Assume (i). 
By the definition of $\SCH_{a}$ there is ${^*\D}\succ \D$ which realises $\stp_a,$ say by $b.$ By
Lemma~\ref{Claim 6},  ${^*\D}\in \E^K_d\smin \E^K_{d-1},$ so ${^*\D}\in \E^K_d\smin \E^K_{d-1}.$

It follows that $\D\in \E^K_d,$ since $\delta^K(X/\ker)\ge -d$ for all $X\subs {^*\V}.$ 

It remains to see that $\D\notin \E^K_{d-1}.$ Indeed, if it were in $\E^K_{d-1},$ we would have $a'\le \D$ with
$\delta^K(a')=-d',$ $d'\le d-1,$ and by the analysis in \ref{4.3} arrive at the fact that 
$\D$ and  ${^*\D}$ realise $\stp_{a'},$ hence using again  \ref{Claim 6}, ${^*\D}\in \E^K_{d'},$ a contradiction. 
This proves (ii).

Now, conversely, assume (ii). 

We  claim that ${^*\D}\in \E^K_d.$ Indeed, since $\D\in \E^K_d$ by \ref{4.3}
we find a $b$  with $\delta^K(b/\ker)=-d,$ so $b\le \D$ and $\D\models \SCH_b.$ So ${^*\D}\models \SCH_b$ and as shown
in the first part of the proof, it follows ${^*\D}\in \E^K_d.$
   
By assumptions, up to isomorphism, $a$ is in ${^*\D}.$ Since $\delta^K(a)=-d,$ we have $a\le  {^*\D}.$
It follows ${^*\D}\models \SCH_{a},$ so  ${\D}\models \SCH_{a}$ and (i) proved.

The second statement of the theorem now follows from Proposition~\ref{2.9}.
\qed


\epk

\bpk  \label{bzeroes} {\bf Theorem.} Let $L\subs \V^n$ be a $K$-linear subspace of dimension $\l$ and
$W\subs (\F^\times)^{n+m}$ be a variety over $\Q.$  Define the set of special parameters   
$$S_{\l}(W)=\{ s\in \F^m: \dim W(s)<n-\l\}.$$

Assume $\D\models {\rm PK+}\SCH_{a}.$ Then
there are a number $N$ and $n-1$-dimensional $\Q$-linear subspaces $M_1,\ldots,M_N\subsetneq \V^{n}$ depending on $L$ and $W$  
 such that  for every $s\in S_{\l}(W)$ for some  $a_1,\ldots a_N\in \V^n$ 
$$ L\cap \ln W(s) \subs \bigcup_{i\le N}(M_i+a_i+\ker^n).$$
Moreover, assuming that $L$ is not contained in a proper $\Q$-linear subspace and
that $\D$ has standard kernel, that is $\ker_{\F}=\omega\Z$ for some transcendental $\omega,$
 we have   for every $s\in S_{\l}(W)$  some  $a_1,\ldots a_N\in \V^n$
 $$ L\cap \ln W(s) \subs \bigcup_{i\le N}(M_i+a_i).$$   
     
\pf First we consider the case of a single $s\in S_{\l}.$ Choose a finite set $B\subset \V$ so that $s\subs\ex(B)$ and $B\cup \ker
\le \D.$
Let $z=\la z_1,\ldots,z_n\ra\in L$  such that $\ex(z)\in W(s).$ By the choice of $B$ $$\delta^K(z/B\cup \ker)\ge 0.$$ But $\ld_K(z/B\cup \ker)+\trd(\ex z/\ex B)<n.$ It follows, $\ld_{\Q}(z/B\cup \ker)<n.$ In other words, $$m_1z_1+\ldots m_nz_n-b\in 
\ker,$$ for some  $m_1,\ldots,m_n\in \Z,$ not all zero, and $b\in \ssp_{\Q}(B).$ Denote 
$$M=\{ \la x_1,\ldots,x_n\ra\in \V^n: \ m_1x_1+\ldots m_nx_n=0\}.$$
We have proved that 
\be \label{z} z\in L\ \&\ \ex(z)\in W(p)\ \Rightarrow 
                    z\in M+a+\ker^n\ee
 for some $a\in \ssp_{\Q}(B)^n$  and $M.$

\medskip

Claim. For a given $s$ there is finitely many $M$ and $a$ such that (\ref{z}).

Indeed, if not then the type saying that $z\in L\ \&\ \ex(z)\in W(s)$ and $z\notin M+a+\ker^n,$ for $M$ running through all $\Q$-linear subspaces of codimension 1 and
 $a\in \ssp_{\Q}(B)^n$ is consistent. This type would be realised in some ${^*\D}\succ \D$ contradicting (\ref{z}).
 
 The proved claim implies the existence of the bound $N_s$ on the number of cosets $M+a+\ker^n$ satisfying (\ref{z}). We need to show that there is an $N$ that bounds all the $N_s.$ Assuming such a bound does not exist we can find a $s\in  S_{\l}(W)$ in some ${^*\D}\succ \D$ for which no finite bound $N_s$ does exist, contradicting the Claim. This proves the first part of the theorem.
 
 The ``moreover'' statement is proved in \cite{Z2}, Thm 2, page 36. The assumption of the uniform
 Schanuel's conjecture is only used to quote Corollary 2 of that paper, the statement of which coincides with the first part of the present theorem.
 \qed

 {\bf Remark.} $S_{\l}(W)$ is a constructible subset, that is quantifier-free definable in the field language.
 
 \epk 
 \bpk \label{oldt2}
{\bf Corollary.} {\em Suppose $\C^K,$ the structure on  complex numbers, for some $K\subs \C$ satifies the assumptions
\ref{assA} for some $a$ in $\C,$ that is 
$$\C^K\models {\rm PK+}\SCH_{a}.$$  
Let $L$ be a $K$-linear subspace of $\C^n$ which is not contained in any proper $\Q$-linear subspace of $\C^n,$ and let $W$ be an algebraic subvariety of $\C^{n+m}.$

Then  there are $n-1$-dimensional $\Q$-linear subspaces $M_1,\ldots,M_N\subsetneq \C^{n},$ $N=N(L,W),$  
 such that  for every $s\in S_{\l}(W)$ there are  $a_1,\ldots a_N\in \C^n$ with the property that
every irreducible component  of the analytic set $  L\cap \ln W(s)\}$  is a subset of   
  $M_j+a_j$ for some $j\in \{1,\ldots,N\}.$}

\epk

\ssn{Completeness, near model completeness and superstability} 

\bpk \label{ID} {\bf Proposition} {\em An $\omega$-saturated model of ${\rm PK+}\SCH_{a}$ is of infinite dimension}.

{\bf Proof.} Let $\D$ be a saturated model of the axioms. We assume that type $\stp_a$ is realised by $a.$

We need to find for every $n$ an $n$-tuple $c_1,\ldots,c_n$ such that for every $b_1,\ldots,b_m,$
$$\delta^K(c_1,\ldots,c_n,b_1,\ldots,b_m/a)\ge n$$
 equivalently, assuming  $c_1,\ldots,c_n,b_1,\ldots,b_m$ are $\Q$-linearly independent over $a,$
 $$\ld_K(c_1,\ldots,c_n,b_1,\ldots,b_m/a)+$$ $$+\trd(\ex c_1,\ldots,\ex c_n,\ex b_1,\ldots,\ex b_m/\exp a)\ge n+m.$$
 
 It is enough to find 
 for any special pairs $(L_1,V_1), \ldots,(L_{\l},V_{\l})$ in the $n+m$-space, 
 elements $c_1,\ldots,c_n$ such that 
 for every $b_1,\ldots,b_m,$ $\Q$-linearly independent over $\{ c_1,\ldots,c_n\}\cup a,$
 for each $j=1,\ldots,\l,$ either $\la c_1,\ldots,c_n,b_1,\ldots,b_m\ra\notin L_j$ or
 $\la \ex c_1,\ldots,\ex c_n,\ex b_1,\ldots,\ex b_m\ra\notin V_j.$
 
{\bf Claim.}  There is a number $k(j)$ and  proper $\Q$-linear subspaces $N_{i,j}\subset \F^n,$ $i=1,\ldots,k(j)$ such that for
 any  $c_1,\ldots,c_n,$ for
 any  $b_1,\ldots,b_m,$ $\Q$-linearly independent over $\{ c_1,\ldots,c_n\} \cup a,$ if
 $\la c_1,\ldots,c_n,b_1,\ldots,b_m\ra\in L_j$ and $\la \ex c_1,\ldots,\ex c_n,\ex b_1,\ldots,\ex b_m\ra\in V_j$
 then $ \la c_1,\ldots,c_n\ra \in N_{ij}+\ssp a$ for some $i\le k(j).$
 
 Indeed, otherwise, in a saturated model we will have  $c_1,\ldots,c_n,b_1,\ldots,b_m$
 such that $\ld_\Q(c_1,\ldots,c_n,b_1,\ldots,b_m/a)=n+m$
 and
 $$\ld_K(c_1,\ldots,c_n,b_1,\ldots,b_m/a)+$$ $$+\trd(\ex c_1,\ldots,\ex c_n,\ex b_1,\ldots,\ex b_m)\le \dim L_j+\dim V_j
 < n+m$$
 which contradicts the fact
 $\delta^K(c_1,\ldots,c_n,b_1,\ldots,b_m/a)\ge 0$ established in \ref{p1}. Claim proved.
 
 Now, let $c_0$ be a non-zero element of $\ker$ and $k\in K\setminus \Q.$ Set $c_1=kc_0.$ Let $p_2,\ldots,p_n$ be integers which we will define later and let 
 $c_i=p_ic_1,$ for $i=2,\ldots,n.$ 
 
We choose the $p_i$ so that
$  \la c_1,\ldots,c_n\ra\notin \bigcup_{ij}N_{ij}.$

These are as required. \qed

\epk

\bpk
\df The extension of the initial language $\LL_K$ by existential predicates  
$$E_P(\bar x)\equiv \exists  \bar y P(\bar x,\bar y),$$
where $P$ is a quantifier-free formula, 
is denoted $\LL_K^E.$ 

We assume throughout that $a$ is a tuple in some $\D$ and  $\delta^K(a)=-d$ \epk
\bpk \lb{Claim 7} {\bf Lemma.} {\em Assuming  $\D_1\subs \D$ as $\LL_K^E$-structures and  $\D\models {\rm PK+}\SCH_{a},$ we
have 
 $\D_1\models {\rm PK+}\SCH_{a}$ and $\D_1\le \D.$}
\smallskip

\pf
 $\D_1\in \E^K_d$ for every $\LL_K$-substructure of $\D,$
since facts of the form $\delta^K(X)=m$ are fixed by quantifier-free types.

To see that $\D_1\le \D$ it is enough to show that for a finite $B$
$$B\le \D_1 \Rightarrow  B\le \D.$$
This follows from Proposition~\ref{p1} if we take into account that
$\stp_B$ is  $\LL_K^E$-quantifier-free.

It remains to see that an elementary extension ${^*\D_1}$ of $\D_1$ contains a copy of $A.$ 
This is immediate by the fact that the condition on consistency of $\qftp_A$ is given by existential $\LL_K$-formulas, so  
that is by  $\LL_K^E$-quantifier-free ones.
\qed

\epk

\bpk \lb{c} {\bf Lemma.} {\em Assume $\D_1,\D_2\in \E^K_d$ and $\D_1\models {\rm EC}.$ 
 Suppose $\D_2\le \D_1.$ Then $\D_2\subs \D_1$ in the language $\LL_K^E.$}
\smallskip

\pf Recall that by Proposition~\ref{c3} $\D_1\in \EC.$
Let $a\subs \D_1$ be finite and suppose 
$\D_2\models \exists  y\ P(a,y),$ where
$ P(x,y)$ is quantifier-free. By the definition of $\EC_d$ we get then
$\D_1\models \exists  y\ P(a,y).$ \qed
\epk 
\bpk \lb{cc} {\bf Corollary.} {\em For $\D_1,\D_2\in \EC_{d}$
$$\D_1\subs \D_2\mbox{ as $\LL_K^E$-structures \ \ iff \ \ \ }\D_1\le \D_2.$$}

We then have by Proposition~\ref{p1}.
\epk
\bpk \lb{ax+} {\bf Corollary.}   
$$\D\in \EC_{d} \mbox{ if and only if } \D\models {\rm PK +} \SCH_{a}{\rm +EC}$$
\epk

We say that a (partial) map $\varphi: \D_1\to \D_2$ is an
{\bf $\LL^E_K$-monomorphism}, if it is injective and for any $k$-ary $\LL^E_K$-predicate
$S$ and any $k$-tuple $a$ from the domain of $\varphi$ $$\D_1\models
S(a)\mbox{ \ iff \ }\D_2\models
S(\varphi(a)).$$ 

\bpk \lb{pi} {\bf Lemma.} {\em Let $\D_1$ and $\D_2$ satisfy ${\rm PK} + \SCH_{A}+{\rm EC},$ and $B_1\le 
\D_1,$
$B_2\le \D_2$ such that there is an $\LL_K$-monomorphism
$$\varphi:B_1\to B_2.$$
Let $\D_{B_1}$ and $\D_{B_2}$ be the expansions of $\D_1,$
$\D_2$ by  constants naming elements of $B_1$ and $B_2$ in 
correspondence with $\varphi.$ Then $$\D_{B_1}\equiv \D_{B_2}.$$ }
\smallskip

\pf We prove that given $\omega$-saturated elementary extensions ${^*\D_1}$ of
$\D_1$ and ${^*\D_2}$ of $\D_2,$ given finite $C\subs {^*\D_1},$
$c\in {^*\D_1}$ and a 
$\LL^E_K$-monomorphism $\varphi$ of $B_1\cup C$ into ${^*\D_2}$ one can
extend the monomorphism to $c.$ By symmetry, this yields a winning strategy for
the Ehrenfeucht-Fraisse game, and we are done.

We may assume that  $\varphi$ is the identity and $B_1\cup C=B=\varphi(B).$
It is enough to show that under the assumption for any
 $c\in {^*\D_1}$ we can extend $\varphi$ to some $B'\sups Bc$ as an
$\LL_K$-monomorphism and $B'\le {^*\D_1},$  $\varphi(B')\le {^*\D_2}.$

If $\partial(c/B)=1$ then  define $B'=Bc$ and $\varphi(c)$ to be any element
from ${^*\D_2}$
which is not in the $\partial$-closure of $A$ in ${^*\D_2}$ (use \ref{ID}).
Then $B'$ and $\varphi(B')$ are as required.

If $\partial(c/B)=0$ then extend $c$ to a  finite string $\bar c$ from
 ${^*\D_1}$
 so that $\delta^K(\bar c/B)=0.$  The quantifier free type of $\bar c$
over $B$ is consistent with $\D_1,$ by \ref{c}, and so is realised in ${^*\D_1},$ 
by $\bar b$ say.
 Since  $\delta^K(\bar b/B)=0,$ we have
$A\bar b\le {^*\D_2}$. So, we can define $B'=B\bar c$ and $\varphi(\bar c)=\bar b.$ \qed
\epk

\bpk \lb{at} {\bf Lemma.} {\em Let $\D_1,\D_2$ be  $\omega$-saturated models of ${\rm PK} + \SCH_{A}+{\rm EC},$ $B_1, B_2$ finite subsets of   $\D_1,$ $\D_2,$ 
correspondingly, and $\varphi:B_1\to B_2$ is a  $\LL_K^E$-monomorphism. 
Then, there exists a finite subset $\tilde B_1$ such that
 $\tilde B_1\le \D_1$
and
 $\varphi$ can be extended to  $\tilde B_1$ in such a way that
$$\varphi(\tilde B_1)=\tilde B_2\le \D_2.$$}
\smallskip

\pf Let $b_1$ be a string of all elements of $B_1$ and $c_1$ a tuple in $\D_1$
such that
$\delta^K(b_1\fr c_1)=\partial( b_1).$ It follows $b_1\fr c_1\le \D_1.$ Set $m=\partial( b_1).$

Let $q^0(x\fr y)$ be the $\LL_K$-quantifier-free type of $b_1\fr c_1.$
Let $ b_2$ be a string in $\D_2$ which corresponds to $ b_1.$ Then the $\LL_K^E$-monomorphism 
 guarantees that the type $q^0(b_2\fr y)$ over $b_2$
is consistent and thus there is a $ c_2$ in $\D_2$ realising the type, in particular 
$\partial(b_2)\le \delta^K( b_1\bar c_1)=m=\partial( b_1).$
By symmetry $\partial( b_2)= m=\partial( b_1).$
 Since $\delta^K( b_2 \fr c_2)=\partial(b_2),$ we have  
$b_2\fr c_2\le \D_2.$
Now Lemma~\ref{c} says that $b_2\fr c_2$ is of the same
$\LL_K^E$-quantifier-free
type as $b_1\fr c_1.$ \qed
\epk

\bpk \lb{qe} {\bf Main Theorem.} {\em Given $\D\in \EC_d,$ let finite  $A\le \D.$ Then the following hold: 

(i) The axioms \ ${\rm PK} + \SCH_{A}+{\rm EC}$  
determine the complete theory ${\rm Th}(\D)$ of $\D.$ 

(ii) The theory ${\rm Th}(\D)$
has quantifier elimination in  language $\LL^E_K.$ 

(iii) ${\rm Th}(\D)$ is superstable. 

(iv) The group structure on 
$\ker$ is stably embedded in $\D$, that is no new relations are induced (using parameters) 
on $\ker$ from $\D.$
}
\smallskip

\pf (i) and (ii) It follows from Lemmas~\ref{pi} (with $B_1\cong A\cong B_2$) and \ref{at}  that the theory is complete and
submodel complete. The latter implies  elimination of
quantifiers (see e.g. Theorem 13.1 of \cite{S}).   

(iii) To prove superstability consider  $\D\in \EC_A$ of cardinality $\lambda.$ 
We want to establish the
cardinality of the set $S(\D)$ of complete 1-types over $\D.$ Let
${^*\D}$ be an elementary extension of $\D$ which realises all $n$-types over
$\D$ for all  $n.$  Let $S^{\#}(\D)$ the set of all complete $n$-types
over $\D$ which are realised in ${^*\D}$ by $n$-tuples
$\bar b=\la b_1,\dots,b_n\ra$ such that
$\delta^K(\bar b/\D)=\partial(b_1/\D).$
It follows   that $\card S(\D)\le \card S^{\#}(\D).$

From general properties of $\le$  we get $\F\bar b\le {^*\D},$
and by Lemma~\ref{c} the $\LL_K^E$-quantifier-free type of $\bar b$ over 
$\D$
is determined by the $\LL_K$-quantifier-free type of that.
 Thus 
$\card S(\D)$ is less or equal to the cardinality of $QFS(\D),$ the set of
all $\LL_K$-quantifier-free complete types over $\D.$ 

We claim that
$\card QFS(\D)\le \lambda + 2^{\omega}.$ 
Indeed, each quantifier-free $\LL_K$-type of $\bar b$ over $\D$
 is uniquely determined by the minimal $K$-affine subspace $L$
over $\D$ containing $\bar b$ and, for each $\l\in \N,$ the minimal
algebraic variety $ W^{\frac{1}{\l}}$ containing $\ex(\frac{\bar b}{\l}).$
Notice that, once $W= W^1$ is known, 
for each $\l$ there is at most $\l^n$ choices of  $W^{\frac{1}{\l}}$
($n=|\bar b|$), all conjugated by torsion elements of $(\F^{\times})^n$
of order $\l.$ This branches into at most $2^{\omega}$ types for each of $\lambda$-many varieties $W.$

(iv) Consider again a saturated model $\D$ of the theory. Let $C\le \D$ be an arbitrary finite self-sufficient set  and let $B=\ker\cup C.$ Clearly, $B\le \D.$ 
We claim first that for every finite tuple $\bar b$ in $B$
the complete $\LL_K$-type of $C\cup\bar b$ is determined by the quantifier-free $\LL_K$-type of the tuple. This is again a direct consequence of $C\cup \bar b\le \D,$ by Lemma~\ref{pi}.   Now, since any type of a tuple in the definable $B$ is equivalent to a $\LL_K$-quantifier-free type, any definable subset of $B^n$ is quantifier-free definable, by compactness. We deduce that any $C$-definable subset of $\ker^n$ is  $\LL_K(C)$-quantifier-free definable, hence any subset of $\ker^n$ definable with parameters is $\LL_K$-quantifier-free definable.

More specifically, let $\bar b$ be $\Q$-linearly independent over $C.$ We claim that then it is $K$-linearly independent over $C,$ which follows from the assumption that $\delta^K(\bar b/C)\ge 0.$ 

It follows that quantifier-free $\LL_K(C)$-formulas without parameters restricted to $\ker$ are Boolean combinations of formulas of the form
$m_1x_1+\ldots+m_nx_n=k_1c_1+\ldots+k_pc_p,$ for some $m_1,\ldots,m_n\in \Z,$ $k_1,\ldots,k_p\in K,$ and of the form $\frac{x}{m}\in \ker$ (equivalently, $\ex(\frac{x}{m} )=1$). 
The latter
can be equivalently rewritten as $\exists y\in \ker \, x= my.$ This is the standard form for core formulas in the theory of the $\Z$-group $(\ker, +,0).$ Which proves that the subsets of $\ker^n$ definable in $\D$ are the same as ones definable in $(\ker, +,0).$
\qed

{\bf Remark.} We note that $\SCH_{a}$ is a set of $\LL_K$-sentences without parameters. The subscript $a$ indicates that it depends on the $\LL_K$-type of the tuple $a.$ 
The complete theory of $\D$ depends not just on $d$ but also on the type of a tuple $a$ with
$\delta^K(a)=-d$ which is consistent with the theory.

If $d=0$ then $a$ is empty (or equal, say to $0\in \V$),   $\SCH_{0}$ consists of $\forall$-formulas and $\EC_0$ is a complete elementary class.

\epk


\ssn{  Raising to powers in the complex numbers}
\bpk \label{ws} Consider the structure $\C^K$ for $K\subs \C.$  Assume Schanuel's conjecture or, more specifically,
its form derived in \ref{caseC}:  
   $$\C^K\in \E^K_d, \mbox{ and }a\le \C^K\mbox{ for some tuple }a.$$
\epk

\bpk \label{6.2}{\bf Theorem.} {\em  Assume the corollary of Schanuel's conjecture in the form \ref{ws}. 
Suppose $K\subs \R$. Then   
$$\C^K\models {\rm PK+}\SCH_a {\rm + EC}.$$ 
In particular, these axioms define the complete theory of the structure.}

\smallskip

\pf {\rm PK} and $\SCH_a$ are immediate by assumptions. 
 
 It remains to establish EC, the exponential-algebraic closedness.  We recall the EC-axioms as given in \cite{Z1}. 
 
 Let $L,L_1,\ldots L_k$ ($k\ge 0$) be  $K$-linear subspaces of $\V^{n+p},$  $L_i\subset L,$ for all $i,$ $a\in \V^{p},$  and
 let $V$ be an algebraic subvariety of $(\C^\times)^{n+p}.$  
 EC states that assuming that the pair $(L(a),V(\exp a))$ is {\em free} and {\em normal} (see \cite{Z1} for the definition) and for any $i\le k,$ $L_i(a)$ is a proper subspace of $L(a),$  there must exist a point $b\in L(a)\cap \ln V(\exp a)$ such that $b\notin \cup_{i\le k}L_i(a).$  

Recall that $L(a)$ can be represented as $L(0)+r(a),$ for some $r(a)\in \V^{n}$ (see \ref{linalg}), and so, by shifting both $L(a)$ and  $ \ln V(\exp a)$ by $-r(a)$ we may assume without loss of generality that $L(a)=L,$ a $K$-linear subspace  of $\V^{n},$ does not depend on a parameter.

The $K$-linear space $\V$ under assumptions of the theorem is just $\C$ as a $K$-vector space. 
In particular, we can represent $\V^{n}=\C^n=\R^n+i\R^n,$ the decomposition into the real and imaginary $K$-subspaces (recall that $K\subset \R$).
We denote 
$$\mathrm{Re}(L)=L\cap \R^n,$$ 
 the real part of $L.$ Note that since $L$ is defined by $\R$-linear equations, $L=\mathrm{Re}(L)+i\mathrm{Re}(L).$ In particular,
 for any subset $B\subseteq \mathrm{Re}(L)$ we have $\mathrm{Re}(L)+iB\subseteq L.$ 
 
 Now we need the following.
 \epk
 \bpk {\bf Lemma.}(cf. \cite{Z2}, section 6, Lemma 5). {\em Assume  $(L, V(c))$ is  free and normal. 
 
 There is  a real number $R$ and an integer $m,$ both depending only on $L$ and $V$ but not $c,$ satisfying the following: there are $m$ $\Q$-affine hyperplanes 
 $H_1,\ldots,H_m\subset L$  such that for any ball $B\subset \mathrm{Re}(L)$ of radius $R$ which does not intersect $\bigcup_{j\le m} \mathrm{Re}(H_j)$ 
 $$(\mathrm{Re}(L)+iB)\cap \ln V(c)\setminus  \cup_{i\le k}L_i(c)\neq \emptyset.$$}
 
 {\bf Proof.} Normality implies in particular that  $\dim L+\dim V(c)\ge n.$
By intersecting with a generic hyperplane  we may assume without loss of generality that $\dim L+\dim V(c)=n.$

 Let $C(L,V)$ be the set of all the $c$ which satisfy the assumption of the lemma. By \cite{Z2}, Corollary 5, page 41,  there is an $R$ such that for any ball $B\subset L$ of radius 
 $R$ there is a
 dense subset $C_B(L,V)$ of $C(L,V)$ with the property that,  for any $c\in C_B(L,V),$
 $$(\mathrm{Re}(L)+iB)\cap \ln V(c)\neq \emptyset$$

Fix an $R$ as above, fix a positive integer $m$ (the value of which we will specify later) and consider a ball $B^*\subset \mathrm{Re}(L)$ of radius $2^m\cdot R$ such that 
$B^*\cap \bigcup_{i\le k} \mathrm{Re}(L_i(c))=\emptyset.$
Let $c\in C(L,V).$  By the above there is a sequence $c_t\in  C_{B^*}(L,V),$ $t\in \N,$
converging to $c.$ We choose a non-principal ultrafilter $D$ on $\N$ and write more generally
the fact of convergence along the ultrafilter as $c=c_t/D.$ 

By the choice of the $\{ c_t:t\in \N\}$ for every $t$ there exists
$$\xi_t\in  (\mathrm{Re}(L)+iB^*)\cap \ln V(c_t).$$
 Thus the limit $\xi_t/D$ of the sequence is either a point $b$ in $ (\mathrm{Re}(L)+iB^*)\cap \ln V(c),$ in which case (the regular case) we are done, or alternatively,  the limit point $\exp(\xi_t/D)$ is not in $(\C^\times)^n$ (the singular case). This situation is analysed in \cite{Z2}, Lemma 4, page 41-42. By (iii) of that lemma, after appropriate transformation of variables, there is a positive $\l<n,$ an algebraic
 variety $W\subset \C^{\l+p}$ depending  on $(L,V)$ but not on $c,$ such that
 $\dim W(c)+\dim \mathrm{pr}L<\l,$ where  $\mathrm{pr}L$ is the projection of $L$ on the first
 $\l$ coordinates. Moreover, $W$ belongs to a finite list $\rho(L,V)$ of varieties.
 
Note that the assumption that $(L,W(c)$ is free implies that $L$ is not contained in a proper $\Q$-linear subspace of $\C^n.$
 By Corollary~\ref{oldt2} of the present paper, for some $N$ depending only on $L$ and $W,$ in the singular case there are  $\Q$-affine
 hyperplanes $M_1+a_1,\ldots,M_N+a_N$ of $L,$ where $a_1,\ldots,a_N$ may depend on $c,$ such that $$\xi_t/D\in (M_1+a_1)\cup\ldots\cup (M_N+a_N).$$
 By finiteness of $\rho(L,V)$ there are a number $m\ge N$ depending only on $L$ and $V$  
 and $\Q$-affine hyperplanes
 $H_1,\ldots,H_m$ which depend on $L,V$ and $c$ only, such that $\xi_t/D\in \bigcup_{j\le m}H_j$ 
or $ (\mathrm{Re}(L)+iB^*)\cap \ln V(c),$ 
 for any choice of $\xi_t$ as above. 
 
 It remains to find a way to choose $\xi_t$ so that 
 $i\mathrm{Re}(\xi_t/D)$ avoids the $i\mathrm{Re}(H_1),\ldots,i\mathrm{Re}(H_m).$ We observe (see \cite{Z2}, Lemma 6) that  inside  a ball $B^*$ of radius $2^mR$ one can always choose a
 ball  $B$ of radius $R$ which avoids $m$ given hyperplanes 
 $\mathrm{Re}(H_1),\ldots,\mathrm{Re}(H_m).$  Since $B$ is still big enough we may assume that $i\mathrm{Re}(\xi_t)\in iB$ and so only the first case is possible.  This proves the Lemma. \qed 
 
 The theorem follows.\qed
 \epk

\bpk \label{bkw}
We recall the following result by A.Wilkie, J.Kirby and M.Bays. \medskip

 {\bf Theorem.} (\cite{BKW} 1.3) {\em Let $\F_{\ex}$ be any exponential field, 
let $C$ be an ecl-closed subfield of $\F_{\ex}$, and let $\lambda$ be an m-tuple which is
exponentially algebraically independent over $C,$ $K=\Q(\lambda).$ Then for any tuple $z$ from $\F:$
\be \label{wkb}\trd(\exp(z)/C(\lambda)) + \ld_K(z/ \ker) - \ld_{\Q}(z/ \ker) \ge 0.\ee
In particular, this holds for the exponential field $\C_{\exp}$ of complex numbers and $C= {\rm ecl}(\emptyset).$
}

Here, an {\bf exponential field} $\F_{\ex}$ is  $(\F,+,\cdot,\ex)$ a field structure with a homomorphism $\ex: \F\to \F^{\times}.$
An {\bf ecl-closed subfield} is an exponential subfield $C\subs \F$ that is exponentially-algebraically closed inside $\F_{\ex}$
(see \cite{BKW} for details). 
In the exponential field $\C_{\exp}$ the ecl-closure ${\rm ecl}(X)$ of a countable subset $X$ is  countable, 
by Lemma 5.12 of \cite{Z3}. In particular, all but countably many 
complex numbers are exponentially algebraically independent over ${\rm ecl}(\emptyset).$ \epk

\bpk {\bf Corollary} of (\ref{wkb}).

Since $$\trd(\exp(z)/C, \lambda)=\trd(\exp(z),\lambda/C)-\trd(\lambda/C)\le \trd(\exp z/C)$$
we have as a corollary
$$\trd(\exp(z)/C) + \ld_K(z/ \ker) - \ld_{\Q}(z/ \ker) \ge 0,$$
and a weaker version, which is of interest to us here,
$$\delta^K(z/\ker)=\trd(\exp z) + \ld_K(z/ \ker) - \ld_{\Q}(z/ \ker) \ge 0,$$
which amounts to say that
$$\D\in \E_0.$$

\epk


\bpk {\bf Corollary.} {\em Let a finite subset $\lambda\subs \C$  be exponentially-algebraically independent over ${\rm ecl}(\emptyset)$  and 
let $K=\Q(\lambda).$ Then  $\C^K$ 
satisfies ${\rm PK+}\SCH_0,$ where $\SCH_0$ denotes $\SCH_a$ with $a=0.$

In particular, the statement of Corollary~\ref{oldt2} holds for $\C^K.$ }
\epk


\end{document}